\documentclass{amsart}

\usepackage{amsmath,amssymb,amsfonts,amsthm}
\usepackage{tikz}
\usetikzlibrary{arrows,shapes}
\usetikzlibrary{decorations.markings}
\tikzstyle{every picture}=[baseline=0,heightone]

\theoremstyle{theorem}
\newtheorem{theorem}{Theorem}
\newtheorem{lemma}[theorem]{Lemma}
\newtheorem{proposition}[theorem]{Proposition}
\newtheorem{corollary}[theorem]{Corollary}
\newtheorem*{conjecture*}{Conjecture}

\theoremstyle{definition}
\newtheorem{definition}[theorem]{Definition}

\newtheorem{example}[theorem]{Example}
\newtheorem*{example*}{Example}

\def\tensor{\otimes}

\def\qptmx#1#2{\left(\begin{smallmatrix}#1\\ #2\end{smallmatrix}\right)}

\def\sgn{\mathrm{sgn}}

\def\tr{\mathsf{tr}}

\def\bse{\mathbf{\hat e}}
\def\C{\mathbb{C}}

\def\bfu{\mathbf{u}}
\def\bfv{\mathbf{v}}
\def\bfw{\mathbf{w}}

\tikzstyle heightone=[scale=.7,shift={(0,-.3)}]
\tikzstyle heightones=[scale=.8,xscale=.35,shift={(0,.1)}]
\tikzstyle heightoneonehalf=[scale=.9,shift={(0,-.2)}]
\tikzstyle heighttwo=[scale=.9,shift={(0,-.4)}]
\tikzstyle heighttwos=[scale=.5,xscale=.6,shift={(0,-.1)}]
\tikzstyle heightthree=[scale=.6,shift={(0,-.9)}]
\tikzstyle heightthrees=[scale=.4,xscale=.7,shift={(0,-.2)}]

\tikzstyle vertex=[circle,draw,fill=black,inner sep=1pt]
\tikzstyle ciliation=[circle,draw=none,fill=red,inner sep=1pt,semitransparent]
\tikzstyle ciliatednode=[vertex,pin={[pin distance=1mm,pin edge={semitransparent,red},ciliation]#1:{}}]

\tikzstyle matrix=[black,thick,circle,draw=blue!50,top color=blue!20,bottom color=black!10,scale=.8,inner sep=1pt]
\tikzstyle small matrix=[matrix,scale=.7]
\tikzstyle vector=[black,thick,rectangle,draw=gray!50!yellow,top color=yellow!30,bottom color=black!10,scale=.8,inner sep=2pt]
\tikzstyle small vector=[vector,scale=.8]
\tikzstyle plain vector=[rectangle,draw=none,fill=white,scale=.7]

\tikzstyle basiclabel=[draw=none,fill=none,shape=rectangle,inner sep=2pt,scale=.8]
\tikzstyle leftlabel=[basiclabel,anchor=east]
\tikzstyle rightlabel=[basiclabel,anchor=west]
\tikzstyle bottomlabel=[basiclabel,anchor=north]
\tikzstyle toplabel=[basiclabel,anchor=south]

\tikzstyle trivalent=[very thick]

\tikzstyle arrowstyle=[blue,semitransparent,scale=2]

\tikzstyle directed=[postaction={decorate,decoration={markings,
    mark=at position .65 with {\arrow[arrowstyle]{stealth}}}}]
\tikzstyle reverse directed=[postaction={decorate,decoration={markings,
    mark=at position .65 with {\arrowreversed[arrowstyle]{stealth};}}}]

\tikzstyle with matrix=[postaction={decorate,decoration={markings,
    mark=at position .5 with {\node[matrix]{#1};}}}]
\tikzstyle with small matrix=[postaction={decorate,decoration={markings,
    mark=at position .5 with {\node[small matrix]{#1};}}}]

\tikzstyle directed matrix=[postaction={decorate,decoration={markings,
    mark=at position .9 with {\arrow[arrowstyle]{stealth}},
    mark=at position .35 with {\node[matrix]{#1};}}}]
\tikzstyle directed small matrix=[postaction={decorate,decoration={markings,
    mark=at position .9 with {\arrow[arrowstyle]{stealth}},
    mark=at position .35 with {\node[small matrix]{#1};}}}]
\tikzstyle reverse directed matrix=[postaction={decorate,decoration={markings,
    mark=at position .4 with {\arrowreversed[arrowstyle]{stealth};},
    mark=at position .65 with {\node[matrix]{#1};}}}]
\tikzstyle reverse directed small matrix=[postaction={decorate,decoration={markings,
    mark=at position .4 with {\arrowreversed[arrowstyle]{stealth};},
    mark=at position .65 with {\node[small matrix]{#1};}}}]

\tikzstyle dotdotdot=[decorate,decoration={markings,
    mark=at position .3 with{\node{.};},
    mark=at position .5 with {\node{.};},
    mark=at position .7 with {\node{.};}}]

\tikzstyle wavyup=[out=90,in=-90]
\tikzstyle wavydown=[out=-90,in=90]
\tikzstyle cup=[out=-90,in=-90,looseness=1]
\tikzstyle cup2=[out=-90,in=-90,looseness=.75]
\tikzstyle cup3=[out=-90,in=-90,looseness=.5]
\tikzstyle capp=[out=90,in=90,looseness=1]
\tikzstyle cap2=[out=90,in=90,looseness=.75]
\tikzstyle cap3=[out=90,in=90,looseness=.5]

\tikzstyle permutation=[rectangle,fill=black,draw=black]
\tikzstyle symmetrizer=[rectangle,fill=black,draw=black]
\tikzstyle antisymmetrizer=[rectangle,fill=gray!10,draw=black]

\title{On A Diagrammatic Proof of the Cayley-Hamilton Theorem}
\author[Elisha Peterson]{Elisha Peterson \\
Department of Mathematical Sciences, United States Military Academy, West Point, NY 10996-1905; Phone: 845-938-5659; Fax: 845-938-2409, \texttt{elisha.peterson@usma.edu}}

\begin{document}

\begin{abstract}
This note concerns a one-line diagrammatic proof of the Cayley-Hamilton Theorem. We discuss the proof's implications regarding the ``core truth'' of the theorem, and provide a generalization. We review the notation of trace diagrams and exhibit explicit diagrammatic descriptions of the coefficients of the characteristic polynomial, which occur as the $n+1$ ``simplest'' trace diagrams. We close with a discussion of diagrammatic polarization related to the theorem.
\end{abstract}

\maketitle


\section{Introduction}

Given an $n\times n$ matrix $A$, $p(\lambda)=\det(A-\lambda I)$ is a degree $n$ polynomial in $\lambda$ called the \emph{characteristic polynomial} of $A$, whose roots are the eigenvalues of the matrix. The \emph{Cayley-Hamilton theorem} says that a matrix satisfies its own characteristic equation $p(A)=0$. Several proofs of this result have been given, including recent approaches via power series \cite{Ber09} and adjugates \cite{deOl07}.


The purpose of this paper is to explain and generalize the following one-line proof of the Cayley-Hamilton Theorem, which appears in Section 6.5 of \cite{Cvi08}:
\begin{theorem}[Diagrammatic Cayley-Hamilton theorem]\label{t:ch-diagrammatic}
    Let
    $\tikz[heightones,xscale=.75]{
        \foreach\xa in {1,2,5}{\draw(\xa,0)to(\xa,1);}
        \foreach\xa in {.1,.9}{\draw[dotdotdot](2,\xa)to(5,\xa);}
        \draw[antisymmetrizer](.7,.3)rectangle node[basiclabel,scale=.8]{$n+1$}(5.3,.7);}$
    represent the antisymmetrizer on $n+1$ vectors in $\C^n$.
    If $A$ is an $n\times n$ matrix, then
        \begin{equation}\label{eq:ch-diagrammatic}
        \tikz[xscale=.35]{\foreach\xa in {0,1,3}{\draw(\xa,0)to(\xa,1);}
        \draw[antisymmetrizer](-.3,.3)rectangle node[basiclabel]{$n+1$}(3.3,.7);
        \draw(0,-.5)to(0,0)(0,1)to(0,1.5)
            (1,1)node[small matrix]{$A$}to[out=90,in=90,looseness=.5](4.5,1)--(4.5,0)to[out=-90,in=-90,looseness=.5](1,0)
            (3,1)node[small matrix]{$A$}to[out=90,in=90,looseness=.75](4,1)--(4,0)to[out=-90,in=-90,looseness=.75](3,0);
        \node[basiclabel]at(2,0){$\cdots$};}
        =0.
        \end{equation}
\begin{proof}
    Since $n+1$ vectors in $\C^n$ are linearly dependent, it follows that
    $\tikz[heightones,xscale=.75]{
        \foreach\xa in {1,2,5}{\draw(\xa,0)to(\xa,1);}
        \foreach\xa in {.1,.9}{\draw[dotdotdot](2,\xa)to(5,\xa);}
        \draw[antisymmetrizer](.7,.3)rectangle node[basiclabel,scale=.8]{$n+1$}(5.3,.7);}
    =0$.
\end{proof}
\end{theorem}
Later sections will explain the notation used in the above theorem, and particularly why \eqref{eq:ch-diagrammatic} is the matrix polynomial equation given in the Cayley-Hamilton Theorem (Proposition \ref{p:characteristic-coefficients}). The proof makes it clear that, from this point-of-view, the Cayley-Hamilton Theorem is a direct consequence of the dimension of the vector space, which leads to linear dependence. Moreover, one is immediately able to generalize the diagrammatic equation to the following, whose proof is just as quick:
\begin{theorem}\label{t:ch-generalized}
    Let $A_1,A_2,\ldots,A_n$ be $n\times n$ matrices. Then
        \begin{equation}\label{eq:ch-generalized}
        \tikz[xscale=.35]{\foreach\xa in {0,1,3}{\draw(\xa,0)to(\xa,1);}
        \draw[antisymmetrizer](-.3,.3)rectangle node[basiclabel]{$n+1$}(3.3,.7);
        \draw(0,-.5)to(0,0)(0,1)to(0,1.5)
            (1,1)node[small matrix]{$A_1$}to[out=90,in=90,looseness=.5](6,1)--(6,0)to[out=-90,in=-90,looseness=.5](1,0)
            (3,1)node[small matrix]{$A_n$}to[out=90,in=90,looseness=.75](5,1)--(5,0)to[out=-90,in=-90,looseness=.75](3,0);
        \node[basiclabel]at(2,0){$\cdots$};}
        =0.
        \end{equation}
\end{theorem}

The diagrammatic exposition here differs substantially from \cite{Cvi08}, where Theorem \ref{t:ch-diagrammatic} is also proven. The additional contribution here is Theorem \ref{t:ch-generalized} and the procedure for generalizing such identities. In addition, trace diagrams provide a more accessible language for describing diagrammatic formulas.

\textit{This paper is organized as follows.} In Section \ref{s:preliminaries}, we review some background in multilinear algebra and signed graph colorings. In Section \ref{s:tracediagrams}, we define trace diagrams and give some of their basic properties. The explanation of the diagrammatic Cayley-Hamilton Theorem is given in Section \ref{s:ch}, and we conclude with a discussion of identity polarization in Section \ref{s:conclusion}.

\section{Preliminaries}\label{s:preliminaries}

\subsection{Multilinear Algebra}\label{ss:multilinear}
This section introduces our notation for multilinear algebra and tensors. A nice introductory treatment of tensors is given in Appendix B of \cite{FH91}.

Let $V$ be a finite-dimensional vector space over a field $F$. 
In what follows, we assume that $V$ has basis $\{\bse_1,\bse_2,\ldots,\bse_n\}$. The space of $k$-tensors $V^{\otimes k}\equiv V\otimes\cdots\otimes V$ is itself a vector space with $n^k$ basis elements of the form
    $$\bse_{\vec\alpha} \equiv \bse_{\alpha_1} \otimes \bse_{\alpha_2} \otimes \cdots \otimes \bse_{\alpha_k},$$
one for each $\vec\alpha=(\alpha_1,\alpha_2,\ldots,\alpha_k) \in N^k$. By convention $V^{\tensor 0}=F$.
Denote by $\mathrm{Fun}(V^{\otimes j},V^{\otimes k})$ the space of \emph{multilinear functions} from $V^{\otimes j}$ to $V^{\otimes k}$. The tensor operation $(f_1,f_2)\mapsto f_1\otimes f_2$ on functions makes the space a \emph{monoidal category}.

\subsection{Signed Graph Coloring}\label{ss:colorings}
This section introduces graph theoretic principles that will be used in defining trace diagrammatic functions. The basic ideas contained here were first given in \cite{Pet06}, and further developed in \cite{MP09}.

Given a vertex $v\in V$ of a graph $G=(V,E)$, denote by $E(v)$ the set of edges adjacent to $v$. We say that two edges in the same $E(v)$ for some $v$ are \emph{adjacent}.

\begin{figure}[htb]
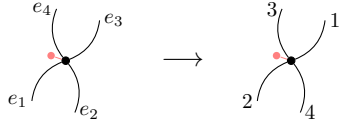

    $$\tikz[shift={(0,.4)}]{
        \node[ciliatednode=170](va)at(0,0){}
            edge[bend right](50:1)edge[bend left](100:1)edge[bend right](230:1)edge[bend left](280:1);
        \node[rightlabel]at(50:1){$e_3$};\node[leftlabel]at(100:1){$e_4$};
        \node[leftlabel]at(230:1){$e_1$};\node[rightlabel]at(280:1){$e_2$};
    }
    \quad\longrightarrow\quad
    \tikz[shift={(0,.4)}]{
        \node[ciliatednode=170](va)at(0,0){}
            edge[bend right](50:1)edge[bend left](100:1)edge[bend right](230:1)edge[bend left](280:1);
        \node[rightlabel]at(50:1){$1$};\node[leftlabel]at(100:1){$3$};
        \node[leftlabel]at(230:1){$2$};\node[rightlabel]at(280:1){$4$};
    }$$
    \caption{By convention, a ciliation on a vertex induces a counter-clockwise ordering, shown as $(e_1,e_2,e_3,e_4)$ at left. The coloring shown at right induces the permutation $\tbinom{2\:4\:1\:3}{1\:2\:3\:4}$ at the vertex.}\label{f:ciliation}
\end{figure}

\begin{definition}
    A \emph{ciliated graph} $G=(V,E,\sigma_*)$ is a graph $(V,E)$ together with an ordering $\sigma_v:\{1,2,\ldots,\deg(v)\}\to E(v)$ of edges at each vertex $v\in V$.
\end{definition}
By convention, when such graphs are drawn in the plane, the ordering is specified by enumerating edges in a counter-clockwise fashion from a ciliation, as shown in Figure \ref{f:ciliation}. Ciliated graphs are sometimes also called \emph{fat graphs}.

\begin{definition}
    Given the set $N=\{1,2,\ldots,n\}$, an \emph{$n$-edge coloring} of a graph $G=(V,E)$ is a map $\kappa:E\to N$ such that no two adjacent edges have the same label. We denote the set of all $n$-edge colorings of a graph $G$ by $\mathsf{col}(G)$ when $n$ is understood.
\end{definition}
In graph theory, edge colorings are sometimes called \emph{Tait colorings}.

Edge colorings induce permutations at the vertices of ciliated graphs. Given an edge coloring $\kappa$ and a degree-$n$ vertex $v$, there is a well-defined permutation $\pi_\kappa(v)\in S_n$ at each interior vertex $v\in V_n$ defined by
    $$\pi_\kappa(v): i \mapsto \kappa(\sigma_v(i)).$$
In other words, 1 is taken to the label on the first edge adjacent to the vertex, 2 is taken to the label on the second edge, and so on, as in Figure \ref{f:ciliation}.

\begin{definition}
    Given an admissible coloring $\kappa$ of a ciliated graph $G=(V,E,\sigma_*)$, the \emph{signature} $\sgn_\kappa(G)$ is the product of permutation signatures on the degree-$n$ vertices:
        $$\sgn_\kappa(G) = \prod_{v\in V_n} \sgn(\pi_\kappa(v)),$$
    where $\sgn(\pi_\kappa(v))$ is the signature of the permutation $\pi_\kappa(v)$.
%
\end{definition}

\begin{definition}
    A \emph{pre-coloring} of a graph $G=(E,V)$ is a coloring $\check\kappa:\check E\to N$ of a subset $\check E\subset E$ of the edges of $G$. A \emph{leaf coloring} is a pre-coloring of the edges adjacent to the degree-1 vertices.

    Given edge sets $\check E_1\subset\check E_2$, if the colorings $\check\kappa_1:\check E_1\to N$ and $\check\kappa_2:\check E_2\to N$ agree on $\check E_1$, we say that $\check\kappa_2$ \emph{extends} $\check\kappa_1$ and write $\check\kappa_2\succ\check\kappa_1$.
\end{definition}

\begin{example*}
    For $n=3$, the pre-coloring
		\tikz[scale=.75]{
            \node[ciliatednode=185](va)at(0,.5){}edge[bend left](-.5,1)edge[bend right](.5,1);
            \node[ciliatednode=175](vb)at(0,0){}edge(va)
                edge[bend right]node[leftlabel,pos=1]{$1$}(-.5,-.5)
                edge[bend left]node[rightlabel,pos=1]{$2$}(.5,-.5);
		}
    has two colorings
    
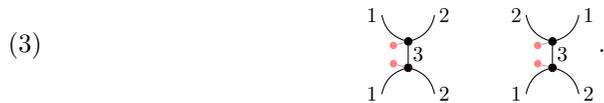
\begin{equation}\label{eq:coloring-ex}
		\tikz{
            \node[ciliatednode=185](va)at(0,.5){}
                edge[bend left]node[leftlabel,pos=1]{$1$}(-.5,1)
                edge[bend right]node[rightlabel,pos=1]{$2$}(.5,1);
            \node[ciliatednode=175](vb)at(0,0){}edge node[rightlabel]{$3$}(va)
                edge[bend right]node[leftlabel,pos=1]{$1$}(-.5,-.5)
                edge[bend left]node[rightlabel,pos=1]{$2$}(.5,-.5);
		}
        \qquad
        \tikz{
            \node[ciliatednode=185](va)at(0,.5){}
                edge[bend left]node[leftlabel,pos=1]{$2$}(-.5,1)
                edge[bend right]node[rightlabel,pos=1]{$1$}(.5,1);
            \node[ciliatednode=175](vb)at(0,0){}edge node[rightlabel]{$3$}(va)
                edge[bend right]node[leftlabel,pos=1]{$1$}(-.5,-.5)
                edge[bend left]node[rightlabel,pos=1]{$2$}(.5,-.5);
		}.
    \end{equation}
    In the first case, the permutations are $\qptmx{1&2&3}{1&2&3}$ and $\qptmx{3&2&1}{1&2&3}$, so that the signature is $\sgn\qptmx{1&2&3}{1&2&3}\cdot\sgn\qptmx{3&2&1}{1&2&3}=-1$. In the second case, the permutations are $\qptmx{1&2&3}{1&2&3}$ and $\qptmx{1&2&3}{3&1&2}$, indicating a positive signature.
\end{example*}

\section{Trace Diagrams}\label{s:tracediagrams}

Penrose was probably the first to describe how tensor algebra may be performed diagrammatically \cite{Pen71}. In his framework, lines in a graph represent elements of the vector field $F$, and nodes represent multilinear functions. Trace diagrams are a generalization of Penrose's tensor diagrams, in which the edges are directed and may be labeled by matrices and the nodes represent the determinant form. The closest concept in graph theory is the \emph{gain graph}, in which edges of a graph are marked by group elements.

\subsection{Definitions}

\begin{definition}\label{d:tracediagrams}
    A \emph{trace diagram} is a directed ciliated graph $\mathcal{D}=(V_1\sqcup V_n,E,\sigma_*)$ together with a marking of edges by linear transformations in $\mathrm{Fun}(V,V)$. Vertices have either degree 1 (in $V_1$) or degree $n$ (in $V_n$). The diagram is \emph{closed} if $V_1$ is empty. A \emph{framed trace diagram} is a diagram together with a partition $E_1=E_I\sqcup E_O$ of the leaf edges adjacent to $V_1$ into ordered \emph{inputs} $E_I$ and \emph{outputs} $E_O$.
\end{definition}
By convention, framed trace diagrams are drawn with inputs at the bottom of the diagram and outputs at the top. Both are assumed to be ordered left to right.

We also permit multiple markings on the same edge, with the understanding that
    $$
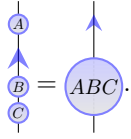
\tikz{
        \draw(0,-.5)to[directed]node[small matrix,pos=.15]{$C$}node[small matrix,pos=.35]{$B$}
            node[small matrix,pos=.825]{$A$}(0,2);}
    =\tikz{\draw(0,-.5)to[directed matrix={$ABC$}](0,2);}.$$
If there are no markings, then $\mathcal{D}$ is a \emph{determinant diagram}, and the orientation of edges may be omitted.

Trace diagrams require an expanded definition of coloring:
\begin{definition}\label{d:tracediagram-coloring}
    A \emph{coloring} of a trace diagram $\mathcal{D}$ is a map $\kappa:E\to N\times N$ labeling the head and tail of each edge by $\kappa(e)_1$ and $\kappa(e)_2$, respectively, in such a way that all labels near an $n$-vertex are different, and unmarked edges $e$ have $\kappa(e)_1=\kappa(e)_2$.
    $$
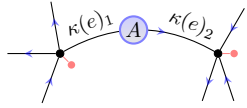
\tikz{
        \node[ciliatednode=-45]at(0,0){}
            edge[directed](100:1)edge[directed](180:1)edge[reverse directed](250:1)
            edge[directed,bend left]node[matrix,pos=.45]{$A$}
            node[toplabel,pos=.2,sloped]{$\kappa(e)_1$}node[toplabel,pos=.8,sloped]{$\kappa(e)_2$}(3,0);
        \node[ciliatednode=0]at(3,0){}
            edge[reverse directed](3.5,.86)edge[directed](2.7,-.9)edge[reverse directed](3.3,-.9);
    }$$
    Suppose further that $A_e$ represents the matrix marking on the edge $e\in E$.
    Then the \emph{coefficient} $\psi_\kappa(\mathcal{D})$ of the coloring is
        $$\psi_\kappa(\mathcal{D})\equiv\prod_{e\in E}(A_e)_{\kappa(e)_2\kappa(e)_1},$$
    where $(A)_{ij}$ represents the $i,j$-matrix entry.
\end{definition}

Thus, the coloring ``picks out'' the entry in the column corresponding to the incoming edge and the row corresponding to the outgoing edge:
    \begin{equation}\label{eq:diagram-matrix-entry}
    \tikz{\draw(0,0)node[rightlabel]{$j$}to[directed small matrix={$A$}](0,1.4)node[rightlabel]{$i$};}
    \quad\leftrightarrow\quad
    (A)_{ij}=\langle\bse_i,A\bse_j\rangle.
    \end{equation}

Pre-colorings $\check\kappa_{\vec\alpha}:(E_I)_i\mapsto\alpha_i$ on the input edges are in a one-to-one correspondence with basis
elements $\bse_{\vec\alpha}\in V^{\tensor|E_I|}$; we will use $\vec\alpha$ as shorthand for the pre-coloring. Together with an output edge coloring $\vec\beta\in V^{\tensor|E_O|}$, the combination $\vec\alpha\cup\vec\beta$ is a leaf coloring.

\begin{definition}\label{d:diagram-function}
    Given a framed trace diagram $\mathcal{D}$, define the \emph{weight} $\chi_{\vec\gamma}(\mathcal{D})$ of a leaf coloring $\gamma$ by
    \begin{equation}
        \chi_{\vec\gamma}(\mathcal{D})=\sum_{\kappa\succ\vec\gamma} \sgn_{\kappa}(\mathcal{D})\psi_{\kappa}(\mathcal{D}).
    \end{equation}
    Define the \emph{trace diagram function} of $\mathcal{D}$ by linear extension of
    \begin{equation}\label{eq:diagram-function}
        f_{\mathcal{D}}:
        \bse_{\vec\alpha}
        \mapsto \sum_{\vec\beta\in N^{|E_O|}}
            \chi_{\vec\alpha\cup\vec\beta}(\mathcal{D}) \bse_{\vec\beta}.
    \end{equation}
    If the diagram is closed, we define its \emph{value} to be
        $$\chi(\mathcal{D})=\sum_{\kappa\in\mathsf{col}(\mathcal{D})}
        \sgn_{\kappa}(\mathcal{D})\psi_{\kappa}(\mathcal{D}).$$
\end{definition}


We will sometimes abuse notation by using the diagram $\mathcal{D}$ interchangeably with $f_\mathcal{D}\in\mathrm{Fun}(V^{\otimes|E_I|},V^{\otimes|E_O|})$. As discussed in the next section, this is permissible because the operation preserves the monoidal structure. We also write formal linear sums of diagrams to indicate the corresponding sums of functions.

\begin{example*}
    The diagram $\tikz{\draw(0,0)to[wavyup](.1,.75);}$ has no vertices, so the signature is trivially +1, and the input and output must have the same label. This implies that $\tikz{\draw(0,0)to[wavyup](.1,.75);}:\bse_i\mapsto\bse_i,$ verifying that the diagram is the identity.

    The two colorings of \eqref{eq:coloring-ex} describe the action of the underlying diagram on the input $\bse_1\tensor\bse_2$:
		$$\tikz{
            \node[ciliatednode=185](va)at(0,.5){}edge[bend left](-.5,1)edge[bend right](.5,1);
            \node[ciliatednode=175](vb)at(0,0){}edge(va)edge[bend right](-.5,-.5)edge[bend left](.5,-.5);
		}:\bse_1\tensor\bse_2\mapsto(\bse_2\tensor\bse_1-\bse_1\tensor\bse_2).
		$$
    The input $\bse_1\tensor\bse_2$ corresponds to the pre-coloring at the bottom of the diagram, while each term in the output corresponds to a coloring extension whose coefficient is coloring's signature.
\end{example*}


\begin{example*}
    The simplest closed trace diagram with a matrix has edge colorings
        \tikz{\draw[directed](0,.5)circle(.5);\node[small matrix]at(.5,.5){$A$};
            \node[basiclabel]at(.5,0){$i$};}
    for $i=1,2,\ldots,n$. Thus, the diagram's value is
        \begin{equation}\label{eq:trace-diagram}
        \tikz{\draw[directed](0,.5)circle(.5);\node[small matrix]at(.5,.5){$A$};}
        =a_{11}+\cdots+a_{nn}=\tr(A).
        \end{equation}
\end{example*}

\subsection{Trace Diagram Relations and Monoidal Structure}\label{ss:functor}

Denote by $\mathfrak{D}(I,O)$ the free $F$-module over framed trace diagrams with $I$ inputs and $O$ outputs. One may compose elements of $\mathfrak{D}(I_1,O_1)$ with those of $\mathfrak{D}(O_1,O_2)$ by gluing outputs to inputs. Since inputs are drawn at the bottom of a diagram and outputs at the top, composition involves drawing one diagram above another.

One may also define $\mathcal{D}_1\tensor\mathcal{D}_2$ as the diagram placing $D_2$ to the right of $D_1$, making the space of framed trace diagrams a monoidal category.

\begin{theorem}[Theorem 5.1 in \cite{MP09}]\label{thm:categorical}
    The mapping $\mathcal{D}\mapsto f_\mathcal{D}$ of Definition \ref{d:diagram-function} is a functor of monoidal categories.
%
\end{theorem}

Intuitively, this result means that a trace diagram's function may be understood by breaking the diagram up into little pieces and gluing them back together.

\begin{definition}
    A \emph{trace diagram relation} is a summation $\sum_{\mathcal{D}} c_{\mathcal{D}} \mathcal{D}\in\mathfrak{D}(I,O)$ of framed trace diagrams for which $\sum_{\mathcal{D}} c_{\mathcal{D}} f_{\mathcal{D}} = 0$.
\end{definition}
Given the monoidal structure, one can apply trace diagram relations on small pieces of larger diagrams (called \emph{local} relations). Sometimes diagrammatic structures are defined as free summations over diagrams modulo one or more local relations \cite{BFK96}.


Trace diagram relations can be made more general than multilinear relations by relaxation of the framing. Denote by $\mathfrak{D}(m)$ the free $F$-module over tensor diagrams with $m$ ordered leaves. A leaf partition gives a mapping $\mathfrak{D}(m)\to\mathfrak{D}(I,O)$, where $I+O=m$. A \emph{(general) trace diagram relation} is a summation $\sum_{\mathcal{D}} c_{\mathcal{D}} \mathcal{D}\in\mathfrak{D}(m)$ that restricts under some leaf partition to a framed trace diagram relation.

\begin{theorem}[Theorem 5.4 in \cite{MP09}]\label{thm:relations}
    Every leaf partition of a general trace diagram relation gives a framed trace diagram relation.
\begin{proof}
    By Definition \ref{d:diagram-function}, the weights of a function depend only on the leaf labels, and not on the partition of framing of the diagram. Since the weights are the same, the relations do not depend on the framing.
\end{proof}
\end{theorem}

The fact that diagrammatic relations are \emph{independent} of framing is very powerful. One may sometimes read off several identities of multilinear algebra from the same diagrammatic relation.

\begin{example*}
    Let $\bfu,\bfv,\bfw\in\C^3$. One can show the cross product and inner product to be
	$$\bfu\times\bfv =
    \tikz{\node[vertex]at(0,.6){}edge(0,1)edge[bend right]node[vector,pos=1]{$\bfu$}(-.5,0)edge[bend left]node[vector,pos=1]{$\bfv$}(.5,0);}
	\quad\text{and}\quad
	\bfu\cdot\bfv =
    \tikz{\draw(0,0)node[vector]{$\bfu$}to[out=90,in=90,looseness=2](1,0)node[vector]{$\bfv$};}.
    $$
    One can also show that
    \begin{equation}\label{eq:3binor}
    \tikz[heighttwo]{\node[vertex]at(0,.5){}edge[bend left](.5,0)edge[bend right](-.5,0);
        \node[vertex]at(0,1){}edge[bend right](.5,1.5)edge[bend left](-.5,1.5)edge(0,.5);}
    =\tikz{\draw(0,0)to[bend left,looseness=.1](.6,1);\draw(.6,0)to[bend right,looseness=.1](0,1);}
    -\tikz{\draw(0,0)to[bend left,looseness=.1](0,1);\draw(.4,0)to[bend right,looseness=.1](.4,1);}.
    \end{equation}
    Consequently,
    $$
    \tikz[scale=.6]{
        \node[vertex](node)at(.5,1){}
            edge[bend right]node[vector,pos=1]{$\bf u$}(0,0)
            edge[bend left]node[vector,pos=1]{$\bf v$}(1,0);
        \node[vertex](node2)at(2.5,1){}
            edge[bend right]node[vector,pos=1]{$\bf w$}(2,0)
            edge[bend left]node[vector,pos=1]{$\bf x$}(3,0)
            edge[bend right=90](node);
    }
    =
    \tikz[scale=.6]{
        \draw(0,0)node[vector]{$\bf u$}to[bend left=90,looseness=2](2,0)node[vector]{$\bf w$};
        \draw(1,0)node[vector]{$\bf v$}to[bend left=90,looseness=2](3,0)node[vector]{$\bf x$};
    }
    -\tikz[scale=.6]{
        \draw(0,0)node[vector]{$\bf u$}to[bend left=90,looseness=1.7](3,0)node[vector]{$\bf x$};
        \draw(1,0)node[vector]{$\bf v$}to[bend left=90,looseness=2.5](2,0)node[vector]{$\bf w$};
    },
    $$
    which is the vector identity
	$$(\bfu\times\bfv)\cdot(\bfw\times\mathbf{x})
        =(\bfu\cdot\bfw)(\bfv\cdot\mathbf{x})-(\bfu\cdot\mathbf{x})(\bfv\cdot\bfw).$$
\end{example*}

\subsection{Symmetrization and Anti-symmetrization}

\begin{definition}\label{d:diagram-symmetrizer-2}
    Define the \emph{anti-symmetrization diagram} by
    $$\index{antisymmetrizer}
    \tikz[heightones,xscale=.6]{
        \foreach\xa in {1,2,5}{\draw(\xa,0)to(\xa,1);}
        \foreach\xa in {.1,.9}{\draw[dotdotdot](2,\xa)to(5,\xa);}
        \draw[antisymmetrizer](.7,.3)rectangle node[basiclabel]{$k$}(5.3,.7);}
    :\bse_{\vec\alpha} \longmapsto \sum_{\sigma\in S_k} \sgn(\sigma)\bse_{\sigma(\vec\alpha)}.$$
\end{definition}

Note that if $k>n$, symmetrization maps any element to 0, since there cannot be more than $n$ linearly independent elements of $M$. This is the observation that proves Theorems \ref{t:ch-diagrammatic} and \ref{t:ch-generalized}.

Redundancies in strands adjacent to the same two nodes are captured by the following result:
\begin{proposition}[Lemma 6.6 in \cite{MP09}]\label{p:vertex-connection-multiplicity}
    Let $v_1,v_2\in V_n$ be adjacent vertices in a trace diagram $\mathcal{D}$ with $\check E=E(v_1)\cap E(v_2)$, representing their shared edges, such that all edges $e\in \check E$ are marked by the same matrix. Choose a coloring $\kappa$ of $\mathcal{D}$ that restricts to pre-colorings $\check\kappa$ on $\check E$ and $\check\kappa^c\equiv\kappa|_{\check E^c}$ on $\check E^c$. Then
        $$\chi_{\check\kappa^c}(\mathcal{D}) = |\check E|! \: \chi_\kappa(\mathcal{D}).$$
\end{proposition}
The factorial in this result shows up in the next two propositions.

\begin{proposition}[Proposition 6.4 in \cite{MP09}]\label{p:antisymmetrizer-two-node}\index{antisymmetrizer}
    \begin{equation}
    \tikz[xscale=.3,yscale=1.5,shift={(0,-.1)}]{
        \foreach\xa in {1,2,5}{\draw(\xa,0)to(\xa,1);}
        \foreach\xa in {.1,.9}{\draw[dotdotdot](2,\xa)to(5,\xa);}
        \draw[antisymmetrizer](.7,.3)rectangle node[basiclabel]{$k$}(5.3,.7);}
    =\frac{(-1)^{\lfloor\frac{n}{2}\rfloor}}{(n-k)!}
    \tikz[heighttwo]{
        \node[ciliatednode=180](vb)at(0,1.5){}
            edge[bend left](-.7,2.2)edge[bend left](-.4,2.2)edge[bend right](.7,2.2);
        \node[ciliatednode=180](va)at(0,.5){}
            edge[out=160,in=200,looseness=2](vb)edge[out=20,in=-20,looseness=2](vb)edge[out=20,in=-20,looseness=1.2](vb)
            edge[bend right](-.7,-.2)edge[bend right](-.4,-.2)edge[bend left](.7,-.2);
        \draw[dotdotdot](-.4,1.9)to node[toplabel]{$k$}(.6,1.9);
        \draw[dotdotdot](-.6,.8)to node[toplabel,scale=.9,xscale=.9]{$n-k$}(.4,.8);
        \draw[dotdotdot](-.4,.1)to node[bottomlabel]{$k$}(.6,.1);
    }
    \end{equation}
\end{proposition}

Note in particular the case with $k=n$:
    \begin{equation}
    \tikz[xscale=.3,yscale=1.5,shift={(0,-.1)}]{
        \foreach\xa in {1,2,5}{\draw(\xa,0)to(\xa,1);}
        \foreach\xa in {.1,.9}{\draw[dotdotdot](2,\xa)to(5,\xa);}
        \draw[antisymmetrizer](.7,.3)rectangle node[basiclabel]{$n$}(5.3,.7);}
    =(-1)^{\lfloor\frac{n}{2}\rfloor}
    \tikz[heighttwo]{
        \node[ciliatednode=180](vb)at(0,1.4){}
            edge[bend left](-.7,2.2)edge[bend left](-.4,2.2)edge[bend right](.7,2.2);
        \node[ciliatednode=180](va)at(0,.6){}
            edge[bend right](-.7,-.2)edge[bend right](-.4,-.2)edge[bend left](.7,-.2);
        \draw[dotdotdot](-.4,1.9)to node[toplabel]{$n$}(.6,1.9);
        \draw[dotdotdot](-.4,.1)to node[bottomlabel]{$n$}(.6,.1);
    }
    \end{equation}

\begin{proposition}[Determinant Diagram, Corollary 6.5 in \cite{MP09}]\label{p:determinant-diagram}\index{determinant!diagram}\index{trace diagram!determinant}
\begin{equation}\label{detdiagram}
    \tikz[heightoneonehalf]{
            \node[ciliatednode=170](topnode)at(0,1.5){};
            \node[ciliatednode=-170](bottomnode)at(0,0){};
            \draw[directed small matrix={$A$}](bottomnode)arc(270:90:.75);
            \draw[directed small matrix={$A$}](bottomnode)arc(-90:90:.75);
            \draw[directed small matrix={$A$}](bottomnode)to[out=135,in=-135](topnode);
            \draw[dotdotdot](-.25,.75)--(.65,.75);
        }
    =(-1)^{\lfloor\frac{n}{2}\rfloor} n! \det(A).
\end{equation}
\end{proposition}

\section{The Not-So-Characteristic Equation}\label{s:ch}\index{Cayley-Hamilton theorem|(}\index{characteristic polynomial|(}

We now return to Theorem \ref{t:ch-diagrammatic}. The fundamental question is how this relates to the formula $\det(A-\lambda I)=0$. The following example demonstrates the relationship in one case:

\begin{example}\label{ex:ch3}
    For $2\times 2$ matrices, Theorem \ref{t:ch-diagrammatic} implies
    \begin{align*}
    6\:\:\tikz[xscale=.35]{\foreach\xa in {0,1,2}{\draw(\xa,0)to(\xa,1);}
        \draw[antisymmetrizer](-.3,.3)rectangle(2.3,.7);
        \draw(0,-.5)to(0,0)(0,1)to(0,1.5)
            (1,1)node[small matrix]{$A$}to[out=90,in=90,looseness=.5](3.5,1)--(3.5,0)to[out=-90,in=-90,looseness=.5](1,0)
            (2,1)node[small matrix]{$A$}to[out=90,in=90,looseness=.75](3,1)--(3,0)to[out=-90,in=-90,looseness=.75](2,0);}
    &=\tikz[xscale=.35]{\draw(0,0)to[wavyup](.1,.5)to[wavyup](0,1)(1,0)to(1,1)(2,0)to[wavyup](1.9,.5)to[wavyup](2,1);
        \draw(0,-.5)to(0,0)(0,1)to(0,1.5)
            (1,1)node[small matrix]{$A$}to[out=90,in=90,looseness=.5](3.5,1)--(3.5,0)to[out=-90,in=-90,looseness=.5](1,0)
            (2,1)node[small matrix]{$A$}to[out=90,in=90,looseness=.75](3,1)--(3,0)to[out=-90,in=-90,looseness=.75](2,0);}
    +\tikz[xscale=.35]{\draw(0,0)to[wavyup](1,1)(1,0)to[wavyup](2,1)(2,0)to[wavyup,looseness=.5](0,1);
        \draw(0,-.5)to(0,0)(0,1)to(0,1.5)
            (1,1)node[small matrix]{$A$}to[out=90,in=90,looseness=.5](3.5,1)--(3.5,0)to[out=-90,in=-90,looseness=.5](1,0)
            (2,1)node[small matrix]{$A$}to[out=90,in=90,looseness=.75](3,1)--(3,0)to[out=-90,in=-90,looseness=.75](2,0);}
    +\tikz[xscale=.35]{\draw(0,0)to[wavyup,looseness=.5](2,1)(1,0)to[wavyup](0,1)(2,0)to[wavyup](1,1);
        \draw(0,-.5)to(0,0)(0,1)to(0,1.5)
            (1,1)node[small matrix]{$A$}to[out=90,in=90,looseness=.5](3.5,1)--(3.5,0)to[out=-90,in=-90,looseness=.5](1,0)
            (2,1)node[small matrix]{$A$}to[out=90,in=90,looseness=.75](3,1)--(3,0)to[out=-90,in=-90,looseness=.75](2,0);}
    -\tikz[xscale=.35]{\draw(0,0)to[wavyup](1,1)(1,0)to[wavyup](0,1)(2,0)to[wavyup](1.9,.5)to[wavyup](2,1);
        \draw(0,-.5)to(0,0)(0,1)to(0,1.5)
            (1,1)node[small matrix]{$A$}to[out=90,in=90,looseness=.5](3.5,1)--(3.5,0)to[out=-90,in=-90,looseness=.5](1,0)
            (2,1)node[small matrix]{$A$}to[out=90,in=90,looseness=.75](3,1)--(3,0)to[out=-90,in=-90,looseness=.75](2,0);}
    -\tikz[xscale=.35]{\draw(0,0)to[wavyup,looseness=.5](2,1)(2,0)to[wavyup,looseness=.5](0,1)(1,0)to(1,1);
        \draw(0,-.5)to(0,0)(0,1)to(0,1.5)
            (1,1)node[small matrix]{$A$}to[out=90,in=90,looseness=.5](3.5,1)--(3.5,0)to[out=-90,in=-90,looseness=.5](1,0)
            (2,1)node[small matrix]{$A$}to[out=90,in=90,looseness=.75](3,1)--(3,0)to[out=-90,in=-90,looseness=.75](2,0);}
    -\tikz[xscale=.35]{\draw(0,0)to[wavyup](.1,.5)to[wavyup](0,1)(1,0)to[wavyup](2,1)(2,0)to[wavyup](1,1);
        \draw(0,-.5)to(0,0)(0,1)to(0,1.5)
            (1,1)node[small matrix]{$A$}to[out=90,in=90,looseness=.5](3.5,1)--(3.5,0)to[out=-90,in=-90,looseness=.5](1,0)
            (2,1)node[small matrix]{$A$}to[out=90,in=90,looseness=.75](3,1)--(3,0)to[out=-90,in=-90,looseness=.75](2,0);}
    \\ &= \tr(A)^2 I + A^2 + A^2 - \tr(A) A - \tr(A) A - \tr(A^2) I
    \\ &= 2\left( A^2 - \tr(A) A + \det(A) I \right) = 0,
    \end{align*}
    where the last step uses the fact that $\det(A)=\frac{1}{2}\left(\tr(A)^2-\tr(A^2)\right)$.
\end{example}

As is usually the case in diagrammatic statements, the ``hard part'' is demonstrating the equivalence to the traditional construction. The remainder of this section proves the following:

\begin{proposition}\label{p:characteristic-coefficients}
    When the antisymmetrizer
    $\tikz[heightones,xscale=.75]{
        \foreach\xa in {1,2,5}{\draw(\xa,0)to(\xa,1);}
        \foreach\xa in {.1,.9}{\draw[dotdotdot](2,\xa)to(5,\xa);}
        \draw[antisymmetrizer](.7,.3)rectangle node[basiclabel,scale=.8]{$n+1$}(5.3,.7);}$
    in \eqref{eq:ch-diagrammatic} is expanded, the coefficients of $A^i$ are equal to $n!$ times the coefficients of $\lambda^i$ in the characteristic polynomial $\det(A-\lambda I)=0$.
\end{proposition}

First, how can the coefficients of the characteristic polynomial be described diagrammatically?

\begin{lemma}\label{l:determinant-sum}
    Given $n\times n$ matrices $A,B$, the determinant sum $\det(A+B)$ is expressed diagrammatically as
    \begin{equation}
    \det(A+B)=(-1)^{\lfloor\frac{n}{2}\rfloor}\sum_{i=0}^n\frac{1}{i!(n-i)!}
    \tikz[heightoneonehalf]{
            \node[ciliatednode=170](topnode)at(0,1.5){};
            \node[ciliatednode=-170](bottomnode)at(0,0){};
            \draw[directed small matrix={$A$}](bottomnode)to[out=180,in=-180,looseness=2.5](topnode);
            \draw[directed small matrix={$A$}](bottomnode)to[out=120,in=-120](topnode);
            \draw[directed small matrix={$B$}](bottomnode)to[out=60,in=-60](topnode);
            \draw[directed small matrix={$B$}](bottomnode)to[out=0,in=0,looseness=2.5](topnode);
            \draw[dotdotdot](-1.25,.75)--(-.25,.75)node[pos=.5,toplabel,scale=.7]{$n-i$};
            \draw[dotdotdot](1.25,.75)--(.25,.75)node[pos=.5,toplabel,scale=.8]{$i$};
        }.
    \end{equation}
\begin{proof}
    Proposition \ref{p:determinant-diagram} states that
    \begin{equation}
    \tikz[heightoneonehalf]{
            \node[ciliatednode=170](topnode)at(0,1.5){};
            \node[ciliatednode=-170](bottomnode)at(0,0){};
            \draw[directed small matrix={$A$}](bottomnode)arc(270:90:.75);
            \draw[directed small matrix={$A$}](bottomnode)arc(-90:90:.75);
            \draw[directed small matrix={$A$}](bottomnode)to[out=135,in=-135](topnode);
            \draw[dotdotdot](-.25,.75)--(.65,.75);
        }
    =(-1)^{\lfloor\frac{n}{2}\rfloor} n! \det(A).
    \end{equation}
    Replacing $A$ with $A+B$ in this diagram, one obtains $2^n$ diagrams, each of which has the form
    $\tikz[heightoneonehalf]{
            \node[ciliatednode=170](topnode)at(0,1.5){};
            \node[ciliatednode=-170](bottomnode)at(0,0){};
            \draw[directed small matrix={$?$}](bottomnode)arc(270:90:.75);
            \draw[directed small matrix={$?$}](bottomnode)arc(-90:90:.75);
            \draw[directed small matrix={$?$}](bottomnode)to[out=135,in=-135](topnode);
            \draw[dotdotdot](-.25,.75)--(.65,.75);
        }$
    where $?$ is either $A$ or $B$. By Proposition \ref{p:vertex-connection-multiplicity}, one may rearrange the strands of the diagrams so that the $A$'s and the $B$'s are grouped separately, without changing the value of the diagram. Therefore, all $\tbinom{n}{i}=\frac{n!}{i!(n-i)!}$ diagrams that have $i$ strands labeled by $B$ have the same value, and the result follows.
\end{proof}
\end{lemma}

\begin{corollary}
    In terms of diagrams, the characteristic polynomial is
    \begin{equation}\label{eq:characteristic-coefficients}
        \det(A-\lambda I)
        =\sum_{i=0}^n\left(\frac{(-1)^{i+\lfloor\frac{n}{2}\rfloor}}{i!(n-i)!}
        \tikz[heightoneonehalf]{
            \node[ciliatednode=170](topnode)at(0,1.5){};
            \node[ciliatednode=-170](bottomnode)at(0,0){};
            \draw[directed small matrix={$A$}](bottomnode)to[out=180,in=-180,looseness=2.5](topnode);
            \draw[directed small matrix={$A$}](bottomnode)to[out=120,in=-120](topnode);
            \draw[directed](bottomnode)to[out=60,in=-60](topnode);
            \draw[directed](bottomnode)to[out=0,in=0,looseness=2.5](topnode);
            \draw[dotdotdot](-1.25,.75)--(-.25,.75)node[pos=.5,toplabel,scale=.7]{$n-i$};
            \draw[dotdotdot](1.25,.75)--(.25,.75)node[pos=.5,toplabel,scale=.8]{$i$};
        }
        \right)\lambda^i
        =\sum_{i=0}^n c_i \lambda^i,
    \end{equation}
    where
    \begin{equation}\label{eq:characeteristic-coefficients-2}
        c_i=\frac{(-1)^{i+\lfloor\frac{n}{2}\rfloor}}{i!(n-i)!}
        \tikz[heightoneonehalf]{
            \node[ciliatednode=170](topnode)at(0,1.5){};
            \node[ciliatednode=-170](bottomnode)at(0,0){};
            \draw[directed small matrix={$A$}](bottomnode)to[out=180,in=-180,looseness=2.5](topnode);
            \draw[directed small matrix={$A$}](bottomnode)to[out=120,in=-120](topnode);
            \draw[directed](bottomnode)to[out=60,in=-60](topnode);
            \draw[directed](bottomnode)to[out=0,in=0,looseness=2.5](topnode);
            \draw[dotdotdot](-1.25,.75)--(-.25,.75)node[pos=.5,toplabel,scale=.7]{$n-i$};
            \draw[dotdotdot](1.25,.75)--(.25,.75)node[pos=.5,toplabel,scale=.8]{$i$};
        }.
    \end{equation}
\end{corollary}

This means that the coefficients of the characteristic polynomial are, up to a constant factor, the $n+1$ ``simplest''
diagrams with two nodes. Are these also the coefficients of $A^i$ in \eqref{eq:ch-diagrammatic}? The next lemma provides the combinatorial decomposition of the antisymmetrizer that is required to demonstrate this fact.
\begin{lemma}\label{l:symmetrizer-sum}
    For any $k$ with $0\leq k\leq n$,
    \begin{equation}\label{eq:symmetrizer-sum}
        \tikz[xscale=.35]{\foreach\xa in {0,1,3}{\draw(\xa,0)to(\xa,1);}
            \draw[antisymmetrizer](-.3,.3)rectangle node[basiclabel,scale=.8]{$k+1$}(3.3,.7);
            \draw(0,-.5)to(0,0)(0,1)to(0,1.5)
                (1,1)node[small matrix]{$A$}to[out=90,in=90,looseness=.5](4.5,1)--(4.5,0)to[out=-90,in=-90,looseness=.5](1,0)
                (3,1)node[small matrix]{$A$}to[out=90,in=90,looseness=.75](4,1)--(4,0)to[out=-90,in=-90,looseness=.75](3,0);
            \node[basiclabel]at(2,0){$\cdots$};}
        =\sum_{i=0}^k\frac{(-1)^i k!}{(k-i)!}
        \tikz[xscale=.35]{\foreach\xa in {1,3}{\draw(\xa,0)to(\xa,1);}
            \draw[antisymmetrizer](.7,.3)rectangle node[basiclabel,scale=.8]{$k-i$}(3.3,.7);
            \draw(1,1)node[small matrix]{$A$}to[out=90,in=90,looseness=.5](4.5,1)--(4.5,0)to[out=-90,in=-90,looseness=.5](1,0)
                (3,1)node[small matrix]{$A$}to[out=90,in=90,looseness=.75](4,1)--(4,0)to[out=-90,in=-90,looseness=.75](3,0);
            \node[basiclabel]at(2,0){$\cdots$};}
        \tikz{\draw(0,-.5)to[wavyup,with small matrix={$A^i$}](.1,1.5);}
    \end{equation}
\begin{proof}
    Choose a summand corresponding to a permutation $\sigma\in S_{k+1}$. Write the permutation as $\sigma=\tau\nu$ where $\tau$ is the cycle containing the first element and $\nu$ contains the remaining cycles. Then $|\tau|=i+1$ in the summand, since the left strand passes through $i$ 
    strands, and the summand contributes to the $A^i$ term.

    There are $\frac{k!}{(k-i)!}$ ways to select $\tau$ such that $|\tau|=i+1$. In each case, the remainder of the diagram is closed and all choices of $\nu$ can be consolidated into a single
    $\tikz[heightones,xscale=.75]{
        \foreach\xa in {1,2,5}{\draw(\xa,0)to(\xa,1);}
        \foreach\xa in {.1,.9}{\draw[dotdotdot](2,\xa)to(5,\xa);}
        \draw[antisymmetrizer](.7,.3)rectangle node[basiclabel,scale=.8]{$k-i$}(5.3,.7);}$
    term. The sign of a summand is given by $\sgn(\sigma)=\sgn(\tau)\sgn(\nu)=(-1)^i\sgn(\nu)$, and $\sgn(\nu)$ is also incorporated within the
    $\tikz[heightones,xscale=.75]{
        \foreach\xa in {1,2,5}{\draw(\xa,0)to(\xa,1);}
        \foreach\xa in {.1,.9}{\draw[dotdotdot](2,\xa)to(5,\xa);}
        \draw[antisymmetrizer](.7,.3)rectangle node[basiclabel,scale=.8]{$k-i$}(5.3,.7);}$
    . Hence, the coefficient of $A^i$ is $(-1)^i\frac{k!}{(k-i)!}$ times the closed diagram shown.
\end{proof}
\end{lemma}

\begin{proof}[Proof of Proposition \ref{p:characteristic-coefficients}]
	Letting $k=n$ in Lemma \ref{l:symmetrizer-sum} and applying Proposition \ref{p:antisymmetrizer-two-node} shows that
		\begin{align*}
        \tikz[xscale=.35]{\foreach\xa in {0,1,3}{\draw(\xa,0)to(\xa,1);}
            \draw[antisymmetrizer](-.3,.3)rectangle node[basiclabel,scale=.8]{$n+1$}(3.3,.7);
            \draw(0,-.5)to(0,0)(0,1)to(0,1.5)
                (1,1)node[small matrix]{$A$}to[out=90,in=90,looseness=.5](4.5,1)--(4.5,0)to[out=-90,in=-90,looseness=.5](1,0)
                (3,1)node[small matrix]{$A$}to[out=90,in=90,looseness=.75](4,1)--(4,0)to[out=-90,in=-90,looseness=.75](3,0);
            \node[basiclabel]at(2,0){$\cdots$};}
        &=\sum_{i=0}^n\frac{(-1)^i n!}{(n-i)!}
        \tikz[xscale=.35]{\foreach\xa in {1,3}{\draw(\xa,0)to(\xa,1);}
            \draw[antisymmetrizer](.7,.3)rectangle node[basiclabel,scale=.8]{$n-i$}(3.3,.7);
            \draw(1,1)node[small matrix]{$A$}to[out=90,in=90,looseness=.5](4.5,1)--(4.5,0)to[out=-90,in=-90,looseness=.5](1,0)
                (3,1)node[small matrix]{$A$}to[out=90,in=90,looseness=.75](4,1)--(4,0)to[out=-90,in=-90,looseness=.75](3,0);
            \node[basiclabel]at(2,0){$\cdots$};}
        \tikz{\draw(0,-.5)to[wavyup,with small matrix={$A^i$}](.1,1.5);}\\
		&=\sum_{i=0}^n\frac{(-1)^i(-1)^{\lfloor\frac{n}{2}\rfloor}n!}{(n-i)!i!}
        \tikz[heightoneonehalf]{
            \node[ciliatednode=170](topnode)at(0,1.5){};
            \node[ciliatednode=-170](bottomnode)at(0,0){};
            \draw[directed small matrix={$A$}](bottomnode)to[out=180,in=-180,looseness=2.5](topnode);
            \draw[directed small matrix={$A$}](bottomnode)to[out=120,in=-120](topnode);
            \draw[directed](bottomnode)to[out=60,in=-60](topnode);
            \draw[directed](bottomnode)to[out=0,in=0,looseness=2.5](topnode);
            \draw[dotdotdot](-1.25,.75)--(-.25,.75)node[pos=.5,toplabel,scale=.7]{$n-i$};
            \draw[dotdotdot](1.25,.75)--(.25,.75)node[pos=.5,toplabel,scale=.8]{$i$};
        } A^i\\
        &=n!\sum_{i=0}^n c_i A^i.\qedhere
		\end{align*}
\end{proof}
\index{Cayley-Hamilton theorem|)}\index{characteristic polynomial|)}

\section{Conclusion}\label{s:conclusion}
We have now explained the notation underlying Theorems \ref{t:ch-diagrammatic} and \ref{t:ch-generalized}, and seen that in the first case the diagram translates into what is normally understood as the Cayley-Hamilton theorem. One of the beautiful consequences of the proof is that the coefficients of the characteristic polynomial, sometimes called the fundamental matrix invariants, match up precisely with the $n+1$ ``simplest'' diagrams involving a single $n\times n$ matrix:
    $$
        \tr(A)=
        \tikz{\draw[with small matrix={$A$}](0,.5)circle(.5);}
        \qquad\cdots\qquad
        c_i(A) \varpropto
        \tikz[heightoneonehalf]{
            \node[vertex](topnode)at(0,1.5){};
            \node[vertex](bottomnode)at(0,0){};
            \foreach\xa/\xb in{90/2.2,30/.5}{
                \draw[with small matrix={$A$},pos=.7](bottomnode)to[bend left=\xa,looseness=\xb](topnode);}
            \foreach\xa/\xb in{90/2,30/.5}{
                \draw(bottomnode)to[bend right=\xa,looseness=\xb](topnode);}
            \foreach\xa/\xb in{-.9/-.3,.2/.9}{\draw[dotdotdot](\xa,.75)to(\xb,.75);}
        }
        \qquad\cdots\qquad
        \det(A) \varpropto
        \tikz[heightoneonehalf]{
            \node[vertex](topnode)at(0,1.5){};
            \node[vertex](bottomnode)at(0,0){};
            \draw[with small matrix={$A$}](bottomnode)arc(270:90:.75);
            \draw[with small matrix={$A$}](bottomnode)arc(-90:90:.75);
            \draw[with small matrix={$A$}](bottomnode)to[out=135,in=-135](topnode);
            \draw[dotdotdot](-.25,.75)--(.65,.75);
        }
    $$
The alert reader may notice that if $n$ is even, an additional nontrivial diagram may be constructed from a single node.
\begin{conjecture*}
    For a $2n\times 2n$ skew-symmetric matrix $A$, the \emph{Pfaffian} $\mathrm{Pf}(A)$ corresponds to
    \begin{equation}
       \mathrm{Pf}(A) \varpropto
       \tikz[heighttwo]{
            \node[vertex]at(0,0){};
            \foreach\xa in{1.25,.6,.35}{\draw[with small matrix={$A$}](0,\xa)circle(\xa);}
            \draw[dotdotdot](-.5,1)to(-1.1,1.6);\draw[dotdotdot](.5,1)to(1.1,1.6);\draw[dotdotdot](0,1.2)to(0,2.5);
       }.
    \end{equation}
\end{conjecture*}

Theorem \ref{t:ch-generalized} relates to the technique of \emph{polarization}, which transforms a trace identity with a single matrix $A$ of degree $k$ into an identity with $k$ matrices $A_1,\ldots,A_k$. Formally, the polar version of an identity $\tau(A)=0$ is
    \begin{equation}\label{eq:polarization}
        \frac{1}{k!}\frac{\partial}{\partial\lambda_1}\cdots\frac{\partial}{\partial\lambda_k} \tau(\lambda_1 A_1+\cdots+\lambda_k A_k)=0.
    \end{equation}
Theorem \ref{t:ch-generalized} gives the same result as polarization, as can be inferred by the sum formulae of Section \ref{s:ch}. For example, the simple formula $A^2-\mathrm{tr}(A)A+\det(A)I=0$ corresponds to
    \begin{equation*}
    6\:\:\tikz[xscale=.35]{\foreach\xa in {0,1,2}{\draw(\xa,0)to(\xa,1);}\draw[antisymmetrizer](-.3,.3)rectangle(2.3,.7);
        \draw(0,-.5)--(0,0)(0,1)to(0,1.5)
            (1,1)node[small matrix]{$A_1$}to[out=90,in=90,looseness=.5](3.5,1)--(3.5,0)to[out=-90,in=-90,looseness=.5](1,0)
            (2,1)node[small matrix]{$A_2$}to[out=90,in=90,looseness=.75](3,1)--(3,0)to[out=-90,in=-90,looseness=.75](2,0);}
    =\tikz[xscale=.35]{\draw(0,0)to[wavyup](.1,.5)to[wavyup](0,1)(1,0)to(1,1)(2,0)to[wavyup](1.9,.5)to[wavyup](2,1);
        \draw(0,-.5)--(0,0)(0,1)to(0,1.5)
            (1,1)node[small matrix]{$A_1$}to[out=90,in=90,looseness=.5](3.5,1)--(3.5,0)to[out=-90,in=-90,looseness=.5](1,0)
            (2,1)node[small matrix]{$A_2$}to[out=90,in=90,looseness=.75](3,1)--(3,0)to[out=-90,in=-90,looseness=.75](2,0);}
    +\tikz[xscale=.35]{\draw(0,0)to[wavyup](1,1)(1,0)to[wavyup](2,1)(2,0)to[wavyup,looseness=.5](0,1);
        \draw(0,-.5)--(0,0)(0,1)to(0,1.5)
            (1,1)node[small matrix]{$A_1$}to[out=90,in=90,looseness=.5](3.5,1)--(3.5,0)to[out=-90,in=-90,looseness=.5](1,0)
            (2,1)node[small matrix]{$A_2$}to[out=90,in=90,looseness=.75](3,1)--(3,0)to[out=-90,in=-90,looseness=.75](2,0);}
    +\tikz[xscale=.35]{\draw(0,0)to[wavyup,looseness=.5](2,1)(1,0)to[wavyup](0,1)(2,0)to[wavyup](1,1);
        \draw(0,-.5)--(0,0)(0,1)to(0,1.5)
            (1,1)node[small matrix]{$A_1$}to[out=90,in=90,looseness=.5](3.5,1)--(3.5,0)to[out=-90,in=-90,looseness=.5](1,0)
            (2,1)node[small matrix]{$A_2$}to[out=90,in=90,looseness=.75](3,1)--(3,0)to[out=-90,in=-90,looseness=.75](2,0);}
    -\tikz[xscale=.35]{\draw(0,0)to[wavyup](1,1)(1,0)to[wavyup](0,1)(2,0)to[wavyup](1.9,.5)to[wavyup](2,1);
        \draw(0,-.5)--(0,0)(0,1)to(0,1.5)
            (1,1)node[small matrix]{$A_1$}to[out=90,in=90,looseness=.5](3.5,1)--(3.5,0)to[out=-90,in=-90,looseness=.5](1,0)
            (2,1)node[small matrix]{$A_2$}to[out=90,in=90,looseness=.75](3,1)--(3,0)to[out=-90,in=-90,looseness=.75](2,0);}
    -\tikz[xscale=.35]{\draw(0,0)to[wavyup,looseness=.5](2,1)(2,0)to[wavyup,looseness=.5](0,1)(1,0)to(1,1);
        \draw(0,-.5)--(0,0)(0,1)to(0,1.5)
            (1,1)node[small matrix]{$A_1$}to[out=90,in=90,looseness=.5](3.5,1)--(3.5,0)to[out=-90,in=-90,looseness=.5](1,0)
            (2,1)node[small matrix]{$A_2$}to[out=90,in=90,looseness=.75](3,1)--(3,0)to[out=-90,in=-90,looseness=.75](2,0);}
    -\tikz[xscale=.35]{\draw(0,0)to[wavyup](.1,.5)to[wavyup](0,1)(1,0)to[wavyup](2,1)(2,0)to[wavyup](1,1);
        \draw(0,-.5)--(0,0)(0,1)to(0,1.5)
            (1,1)node[small matrix]{$A_1$}to[out=90,in=90,looseness=.5](3.5,1)--(3.5,0)to[out=-90,in=-90,looseness=.5](1,0)
            (2,1)node[small matrix]{$A_2$}to[out=90,in=90,looseness=.75](3,1)--(3,0)to[out=-90,in=-90,looseness=.75](2,0);}
    =0,
    \end{equation*}
which reads as
    $$\mathrm{tr}(A_1)\mathrm{tr}(A_2)I+A_2A_1+A_1A_2-\mathrm{tr}(A_2)A_1-\mathrm{tr}(A_1)A_2-\mathrm{tr}(A_1A_2)I=0.$$
It is not hard to verify that this is the polarized version of the original identity given by \eqref{eq:polarization}.

The same process can be applied in other scenarios. As a simple example, it is just as easy to prove that
    \begin{equation*}
    6\:\:\tikz[xscale=.35]{\foreach\xa in {0,1,2}{\draw(\xa,0)to(\xa,1);}\draw[antisymmetrizer](-.3,.3)rectangle(2.3,.7);
        \draw(0,-.5)--(0,0)(0,1)node[small matrix]{$A$}to(0,1.5)
            (1,1)node[small matrix]{$B$}to[out=90,in=90,looseness=.5](3.5,1)--(3.5,0)to[out=-90,in=-90,looseness=.5](1,0)
            (2,1)node[small matrix]{$C$}to[out=90,in=90,looseness=.75](3,1)--(3,0)to[out=-90,in=-90,looseness=.75](2,0);}
    =
    \tikz[xscale=.35]{\draw(0,0)to[wavyup](.1,.5)to[wavyup](0,1)(1,0)to(1,1)(2,0)to[wavyup](1.9,.5)to[wavyup](2,1);
        \draw(0,1)node[small matrix]{$A$}to[out=90,in=90,looseness=.5](4,1)--(4,0)to[out=-90,in=-90,looseness=.5](0,0)
            (1,1)node[small matrix]{$B$}to[out=90,in=90,looseness=.5](3.5,1)--(3.5,0)to[out=-90,in=-90,looseness=.5](1,0)
            (2,1)node[small matrix]{$C$}to[out=90,in=90,looseness=.75](3,1)--(3,0)to[out=-90,in=-90,looseness=.75](2,0);}
    +\tikz[xscale=.35]{\draw(0,0)to[wavyup](1,1)(1,0)to[wavyup](2,1)(2,0)to[wavyup,looseness=.5](0,1);
        \draw(0,1)node[small matrix]{$A$}to[out=90,in=90,looseness=.5](4,1)--(4,0)to[out=-90,in=-90,looseness=.5](0,0)
            (1,1)node[small matrix]{$B$}to[out=90,in=90,looseness=.5](3.5,1)--(3.5,0)to[out=-90,in=-90,looseness=.5](1,0)
            (2,1)node[small matrix]{$C$}to[out=90,in=90,looseness=.75](3,1)--(3,0)to[out=-90,in=-90,looseness=.75](2,0);}
    +\tikz[xscale=.35]{\draw(0,0)to[wavyup,looseness=.5](2,1)(1,0)to[wavyup](0,1)(2,0)to[wavyup](1,1);
        \draw(0,1)node[small matrix]{$A$}to[out=90,in=90,looseness=.5](4,1)--(4,0)to[out=-90,in=-90,looseness=.5](0,0)
            (1,1)node[small matrix]{$B$}to[out=90,in=90,looseness=.5](3.5,1)--(3.5,0)to[out=-90,in=-90,looseness=.5](1,0)
            (2,1)node[small matrix]{$C$}to[out=90,in=90,looseness=.75](3,1)--(3,0)to[out=-90,in=-90,looseness=.75](2,0);}
    -\tikz[xscale=.35]{\draw(0,0)to[wavyup](1,1)(1,0)to[wavyup](0,1)(2,0)to[wavyup](1.9,.5)to[wavyup](2,1);
        \draw(0,1)node[small matrix]{$A$}to[out=90,in=90,looseness=.5](4,1)--(4,0)to[out=-90,in=-90,looseness=.5](0,0)
            (1,1)node[small matrix]{$B$}to[out=90,in=90,looseness=.5](3.5,1)--(3.5,0)to[out=-90,in=-90,looseness=.5](1,0)
            (2,1)node[small matrix]{$C$}to[out=90,in=90,looseness=.75](3,1)--(3,0)to[out=-90,in=-90,looseness=.75](2,0);}
    -\tikz[xscale=.35]{\draw(0,0)to[wavyup,looseness=.5](2,1)(2,0)to[wavyup,looseness=.5](0,1)(1,0)to(1,1);
        \draw(0,1)node[small matrix]{$A$}to[out=90,in=90,looseness=.5](4,1)--(4,0)to[out=-90,in=-90,looseness=.5](0,0)
            (1,1)node[small matrix]{$B$}to[out=90,in=90,looseness=.5](3.5,1)--(3.5,0)to[out=-90,in=-90,looseness=.5](1,0)
            (2,1)node[small matrix]{$C$}to[out=90,in=90,looseness=.75](3,1)--(3,0)to[out=-90,in=-90,looseness=.75](2,0);}
    -\tikz[xscale=.35]{\draw(0,0)to[wavyup](.1,.5)to[wavyup](0,1)(1,0)to[wavyup](2,1)(2,0)to[wavyup](1,1);
        \draw(0,1)node[small matrix]{$A$}to[out=90,in=90,looseness=.5](4,1)--(4,0)to[out=-90,in=-90,looseness=.5](0,0)
            (1,1)node[small matrix]{$B$}to[out=90,in=90,looseness=.5](3.5,1)--(3.5,0)to[out=-90,in=-90,looseness=.5](1,0)
            (2,1)node[small matrix]{$C$}to[out=90,in=90,looseness=.75](3,1)--(3,0)to[out=-90,in=-90,looseness=.75](2,0);}
    =0,
    \end{equation*}
Upon taking the trace, one obtains
    \begin{multline*}
        \mathrm{tr}(ABC)+\mathrm{tr}(ACB)\\
        =\mathrm{tr}(AB)\mathrm{tr}(C)+\mathrm{tr}(A)\mathrm{tr}(BC)+\mathrm{tr}(B)\mathrm{tr}(CA)-\mathrm{tr}(A)\mathrm{tr}(B)\mathrm{tr}(C).
    \end{multline*}
This is sometimes called the \emph{Fricke sum relation} and features prominently in the invariant theory of $2\times2$ matrices \cite{Bul97,Dre07,Fri1896,Gol09}.


\bibliographystyle{plain}
\bibliography{characteristic}

\end{document}